\documentclass[a4paper, 
11pt
]{article}

\usepackage{subfig}
\usepackage{pgfplots}
\usepackage{pgfplotstable}
\usepackage{mathtools}
\usepackage{multicol}
\usepackage{comment}
\usepackage{booktabs}
\pgfplotsset{compat=1.5}
\usepackage{amssymb}
\usepackage{url}

\usepackage{pdfsync}
\usepackage{float}
\usepackage{tabularx}
\usepackage{enumerate}
\usepackage{array}
\usepackage{xspace}
\usepgfplotslibrary{groupplots}

\usepackage{authblk}

\newcommand{\x}{x}
\newcommand{\mupar}{\mu}
\makeindex             
\begin{document}

\title{Reduced Order Isogeometric Analysis Approach for PDEs in
  Parametrized Domains}

\author[1]{Fabrizio~Garotta\footnote{fabrizio.garotta01@universitadipavia.it}}
\author[2]{Nicola~Demo\footnote{nicola.demo@sissa.it}}
\author[2]{Marco~Tezzele\footnote{marco.tezzele@sissa.it}}
\author[1]{Massimo~Carraturo\footnote{massimo.carraturo01@universitadipavia.it}}
\author[1]{Alessandro~Reali\footnote{alereali@unipv.it}}
\author[2]{Gianluigi~Rozza\footnote{gianluigi.rozza@sissa.it}}

\affil[1]{Department of Civil Engineering and Architecture
Universit\`a degli Studi di Pavia, Pavia, Italy}
\affil[2]{Mathematics Area, mathLab, SISSA, International School of
  Advanced Studies, via Bonomea 265, I-34136 Trieste, Italy} 

\maketitle

\begin{abstract}
    In this contribution, we coupled the isogeometric analysis to a reduced
    order modelling technique in order to provide a computationally efficient
    solution in parametric domains. In details, we adopt the free-form
    deformation method to obtain the parametric formulation of the domain and
    proper orthogonal decomposition with interpolation for the computational
    reduction of the model. This technique provides a real-time solution for
    any parameter by combining several solutions, in this case computed using
    isogeometric analysis on different geometrical configurations of the
    domain, properly mapped into a reference configuration. We underline that
    this reduced order model requires only the
    full-order solutions, making this approach non-intrusive. We present in
    this work the results of the application of this methodology to a heat
    conduction problem inside a deformable collector pipe.
\end{abstract}

\section{Introduction}
\label{sec:intro}

Nowadays, in the industrial and engineering fields, as well as in the biomedical
sciences, fast and accurate simulations are crucial in several
applications, such as, for example, shape design optimization and
real-time patient specific diagnosis and control. To this end many
reduced order modelling (ROM) techniques have been developed in the
last decade~\cite{schilders2008model,chinesta2016model,quarteroni2014ms,salmoiraghi2016advances,
rozza2018advances}. We cite among others reduced basis
methods~\cite{hesthaven2016certified,rozza2007reduced}, proper
orthogonal decomposition (POD)~\cite{christensen1999evaluation},
proper generalized decomposition~\cite{chinesta2013proper}, and
hierarchical model
reduction~\cite{perotto2010hierarchical,perotto2017higamod,baroli2017hi}. Reduced
order modelling can be integrated to various high fidelity methods
such as finite element~\cite{quarteroni2011certified}, spectral
element, or finite volume
methods~\cite{haasdonk2008reduced,stabile2017pod,stabile2018finite}. We do
mention also recent features of the reduced methods able to provide useful
algorithms for uncertainty quantification as well as data science and better
exploitation of high performance computing~\cite{tezzele2018ecmi, ChenQuarteroniRozza2017}.

Reduced order methods allow a fast and reliable approximation
of parameterized PDEs by constructing small-sized
approximation spaces. Using these spaces for the discretization of the
original problem, it is possible to build a reduced order model
that is a sufficiently accurate approximation of the original full
order problem. The fundamental characteristic that makes the method
functional from an engineering and industrial point of view is that the
offline phase (more expensive), where the actual analysis is carried
out, is performed only once in high performance computing (HPC)
structures and then remains. The online phase 
exploits the calculations already performed and therefore a small
computational power, like the one of laptops or portable devices, is
sufficient. This ensures real-time processing of the problem without
having to access HPC facilities for the analysis of new parameters.

ROM is crucial in industrial simulation-based design optimization
problems in naval and nautical
engineering~\cite{tezzele2018ecmi,demo2018shape}, but also in
biomedical applications for coronary bypass~\cite{Ballarin2017,BallarinJCP} and
carotid occlusions~\cite{tezzele2018combined} for example.

\begin{figure}[!h]
\includegraphics[width=\textwidth]{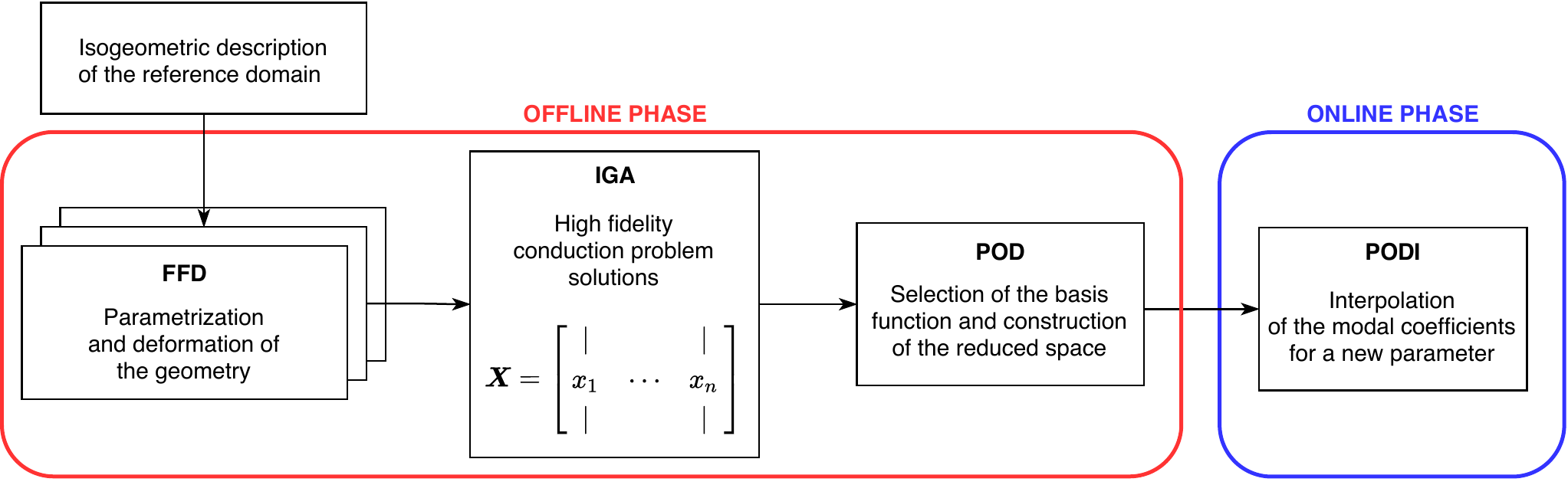}
\caption{Offline-online numerical pipeline. We consider the heat conduction problem
    described by PDEs defined on parametrized geometry. The parametrization of
    the geometry is managed through the FFD which allows the realization of
    different deformation settings. In the offline stage we first solve a
    full order model using IGA to derive the solutions and then, we create a
    ROM applying the POD as space reduction technique. Finally, through the
    PODI, in the online stage we look for the real-time solution of the reduced
    problem for a new parameter.}
    \label{fig:finalscheme}
\end{figure}

The focus of this work is to embed in a ROM framework the isogeometric analysis
(IGA)~\cite{hughes2005isogeometric,cottrell2007studies,cottrell2009isogeometric}
for the simulation of heat diffusion inside a collector pipe. The proposed
approach is  integrated in a numerical pipeline with efficient geometrical
parameterization of the domain through free form deformation
(FFD)~\cite{sederbergparry1986,lassila2010parametric}, an IGA solver as high
fidelity discretization, and POD with interpolation
(PODI)~\cite{bui2003proper,salmoiraghi2018,ripepi2018reduced} for a fast
evaluation of the solution field at untested parameters.
Figure~\ref{fig:finalscheme} depicts the schema of the complete computational
pipeline.

We chose FFD instead of other general purpose geometrical
parameterization techniques such as radial basis functions (RBF)
interpolation~\cite{buhmann2003radial,morris2008cfd,manzoni2012model},
or inverse distance weighting (IDW) interpolation
\cite{shepard1968,witteveenbijl2009,forti2014efficient,BallarinDAmarioPerottoRozza2017},
because of the possibility to use only few parameters to deform the
entire domain of interest. 

The IGA approach allows to integrate classical finite element
analysis (FEA) into conventional industrial CAD tools. To this end IGA
directly employs standard CAD representation bases, e.g., B-splines or
Non-uniform rational B-splines (NURBS), as basis for the analysis. In
this way we can avoid the classical mesh generation and the
consequent geometrical approximation error, obtaining a direct
design-to-analysis simulation, since we are employing the same class
of functions for both the geometry parameterization and the solution
fields approximation.
In this context, the IGA is ideal for
solving elliptic and parabolic PDEs on domains of very general
shape. However, when the objective is to solve the same problem
repeatedly on different domains, the cost of setting up the problem
(meshing, matrix assembly) every time from scratch can be too high. An
optimal solution to this problem is a reduction of the model. 

Previous IGA-ROMs works were developed in the last
years~\cite{manzoni2015reduced,salmoiraghi2016isogeometric,zhu2017isogeometric},
but we underline that the novelty of this work is related to the POD
with interpolation integration into the numerical pipeline, for a
non-intrusive approach. Even if in this work we present a proof of
concept we stress the fact that it can play an important role for
the integration with industrial CAD files being independent from the
IGA full order solver used.

\section{The parametrized heat conduction problem inside a collector pipe}
\label{sec:problem}

The problem of interest we are going to solve throughout this work is
a parametrized heat diffusion problem inside a collector pipe. 

\begin{figure}[h!]
\centering
\includegraphics[width=0.45\textwidth]{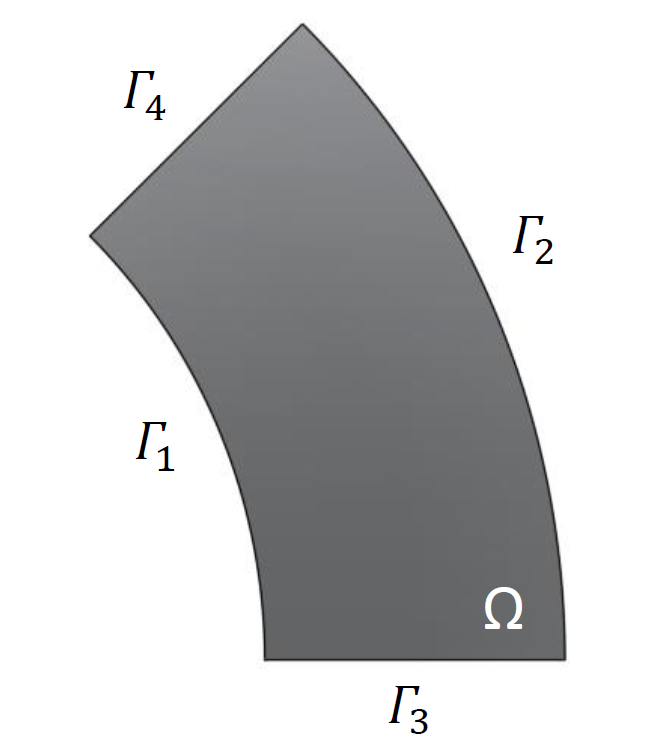}
\caption{Idealized collector pipe representation scheme in
  2D. $\Omega$ represents with internal domain of interest, while
  $\Gamma_{1, \dots, 4}$ indicate the different boundaries. In
  particular $\Gamma_3$ is the inlet, and $\Gamma_4$ is the outlet.}
\label{fig:pipe}
\end{figure}

Let $\Omega \subset \mathbb{R}^2$ be a domain that describes an
idealized collector pipe in 2D, as shown in Figure~\ref{fig:pipe}. We
will refer to $\Omega$ as the reference domain, and for practical
reasons it represents the undeformed geometry. 

We also introduce $\mathbb{D} \subset \mathbb{R}^m$ which is our
parameter space, and for convenience it will be an hypercube. For every $\mupar
\in \mathbb{D}$, which is a vector of geometrical parameters describing
a particular deformation of the domain, we can define a shape morphing map
$\mathcal{M} (\x; \mupar) : \mathbb{R}^2 \to \mathbb{R}^2$. We will
indicate the deformed domain as $\Omega (\mupar) = \mathcal{M}
(\Omega; \mupar)$. We refer to Section~\ref{sec:ffd} for the specific
characterization of such mapping.

The parametrized heat diffusion problem reads: find $u(\mupar)$ such
that
\begin{equation}
  \label{eq:strong}
\begin{cases} \Delta{u} (\mupar) = 0 \quad &\text{in} \;\; \Omega(\mupar) \\
u(\mupar) = 0 \quad &\text{in} \;\; \Gamma_{1,2,4} \\
\nabla u (\mupar) \cdot n = g \quad &\text{in} \;\; \Gamma_3,
\end{cases}
\end{equation}
where $u$ is the temperature distribution inside the domain, and $g$
represents the prescribed heat flux at the inlet. The Dirichlet
boundary conditions describe a perfect insulator with no flux. For
sake of simplicity from now on $g = 1$.

We can introduce the weak formulation of the
problem~\eqref{eq:strong}. We denote with
\[
  V = H^1_{0_{124}}(\Omega) := \left\{v
  \in H^1(\Omega) \; \text{such that} \; v|_{\Gamma_1,\Gamma_2,\Gamma_4} =
  0\right\}
\]
the Sobolev space for the temperature. Multiplying the
first equation of the system by a test function and integrating by
parts we obtain the following problem: given $\mupar \in \mathbb{D}$,
find $u \in V$ such that
\begin{equation}
  \label{eq:weak}
a(u, v; \mupar) = L(v; \mupar) \qquad \forall v \in V,
\end{equation}
where the bilinear form $a(u, v; \mupar)$, and the linear form
$L(v; \mupar)$, are defined as follows

\begin{align}
a(u,v; \mupar) &= \int_{\Omega} \nabla u (\mupar) \; \nabla v \; dV
&&\forall u, v \in V, \\
L(v; \mupar) &= \int_{\Gamma_3} g \; v \; dS &&\forall v \in V.
\end{align}

\section{Isogeometric paradigm for both the geometry and the solution field}
\label{sec:iga}

Usually a CAD representation of the domain is obtained through
B-splines or NURBS, which are able to exactly describe all conic
sections. Here we are going to briefly present both. 

It is possible to derive the B-spline basis functions of order $p$ using Cox-de
Boor's recursion formula~\cite{cox1972numerical,de1972calculating} 
\begin{equation}
  \label{eq:knot3}
N_{i,p}\left(\xi\right)=\frac{\xi-\xi_i}{\xi_{i+p}-\xi_i}N_{i,p-1}
\left(\xi\right)+\frac{\xi_{i+p+1}-\xi}{\xi_{i+p+1}-\xi_{i+1}}N_{i+1,p-1}
\left(\xi\right),
\end{equation}
where
\begin{equation}
  \label{eq:knot2}
N_{i,0}(\xi) =
  \begin{cases}
    1  & \quad \text{if } \quad \xi_i \le \xi < \xi_{i+1},\\
    0  & \quad \text{otherwise,} 
  \end{cases}
\end{equation}
and $\Xi=\{\xi_1,\xi_2,\dots,\xi_{n+p+1}\}$ is the knot vector, a
non-decreasing set of coordinates with $\xi_i
\in \mathbb{R}$, and $n$ is the
number of basis functions which comprise the B-spline.
B-spline curves in $\mathbb{R}^d$ are constructed by taking a linear
combination of B-spline basis functions. The vector valued
coefficients of the basis functions are referred to as \emph{IGA control points}. 
Given $n$ basis functions $\{ N_{i,p} \}_{i=1}^n$ of order
$p$, and the corresponding control points $P_i$, a piecewise polynomial
B-spline curve is given by  
\begin{equation}
  \label{eq:bspline}
C\left(\xi\right)=\sum_{i=1}^{n}{P_iN_{i,p}\left(\xi\right)}.
\end{equation}

IGA control points are points that define the so called control mesh,
which is a mesh made up by the multilinear elements that define and
control the geometry of the problem. It is important to emphasize that
the control mesh does not coincide with the actual geometry of the
physical domain. The control points can be considered as the analog of
the nodal coordinates of the finite element method, with the
difference that, in IGA contest they represent the coefficients of the
basis functions of a B-spline having non-interpolatory nature. 
A generalization of B-splines are the NURBS, which are a rational
version of them and can thus represent exactly any kind of geometry.
This feature of NURBS allows to bypass altogether the computationally
expensive mesh generation and refinement cycle and at the same time to
preserve the exact geometry of the CAD model. 
The key insight of IGA is to use the geometrical map of the NURBS representation as a
basis for the push forward used in the analysis.  
NURBS basis functions of order $p$ are defined through B-spline basis
functions as
\begin{equation}
R_{i,p}(\xi) = \frac{N_{i,p}(\xi)w_i}{W(\xi)}=\frac{N_{i,p}(\xi)w_i}{\sum_{j=1}^{n}N_{j,p}(\xi)w_j},
\end{equation}
where  $w_i$ are the associated weights.
Taking a linear combination of basis functions and control points, we
express a NURBS curve as 
\begin{equation}
C(\xi) = \sum^n_{i=1} P_i R_{i,p}(\xi).
\end{equation}

The isoparametric concept is utilized for both FEA and IGA. However,
the difference between FEA and IGA lies in the bases employed for the
analysis. 
In IGA, the inputs for the calculations come from a CAD model defined
by NURBS curves, which can be used directly for analysis, while in FEA
the finite element mesh is generated starting from an approximation of
the original geometry. 
The mapping from the parametric domain to the physical domain
is then given by 
\begin{equation}
x=\sum_{k=1}^n R_k (\xi)P_k,
\end{equation}
where $R_k(\xi)$ are the NURBS basis functions, $n$ is the number of control
points, $\xi$ the parametric coordinate and $P_k$ is the $k$-th
control point. 
In an isoparametric formulations the displacement field is
approximated by the same shape functions formally: 
\begin{equation}
u=\sum_{k=1}^n R_k (\xi)u_k,
\end{equation}
where $u_k$ is the value of the displacement field at the control
point $P_k$. It is therefore referred to as a control variable or more
generally a degree of freedom. 

\begin{figure} [h!]
\centering
\includegraphics[width=0.75\textwidth]{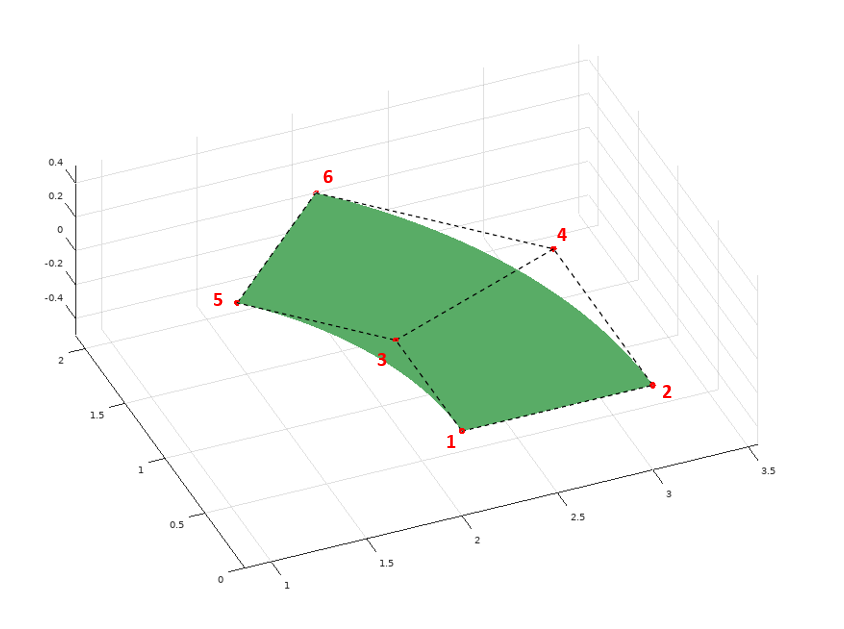}
\caption{Idealized 2D collector pipe and its IGA control points mesh. The six
  control points are indicated with red dots.}
\label{fig:controlp}
\end{figure}

In Figure~\ref{fig:controlp} we present the IGA representation of the
domain $\Omega$ we described in section~\ref{sec:problem}. In red the
six IGA control points defining the NURBS curves. In particular the
knot vectors $\Xi$ and $H$ are defined as follows:  
\begin{align*}
  \Xi &= \left\{0,0,1,1\right\} &\quad k &= 2 \quad p = 1, \\
  H &= \left\{0,0,0,1,1,1\right\} &\quad k &= 3 \quad p = 2,
\end{align*}
where $k$ and $p$ respectively indicate the multiplicity and the degree of the polynomial 
$(k = p + 1)$.

\section{Shape parameterization and deformation through free form deformation}
\label{sec:ffd}

The FFD method has been proposed
in~\cite{sederbergparry1986}.
It was initially used as a tool for computer-assisted geometric design and
animation, nowadays instead it is mostly adopted in academia, industry and
several engineering application fields as morphing technique for complex
geometries thanks to its
features~\cite{demo2018isope,tezzele2018model,tezzele2018dimension,
sieger2015shape}. 
In the FFD procedure, the object to be deformed is embedded into a rectangular
lattice of points, then some of these points are moved to deform the whole
embedded domain. This technique has three main benefits: ({\it i}) with few
parameters --- the displacement of the lattice points --- it is possible to perform global
deformations, ({\it ii}) it allows to preserve continuity also in the surface
derivatives and ({\it iii}) it is completely independent with respect
to the object, so it results applicable also to computational
grids~\cite{lassila2010parametric}.

Initially, FFD maps the original domain $\Omega$ to the reference one using the
affine map $\psi$ defined as $\psi: D \rightarrow [0, 1]^n$, where $D \supset
\Omega$ is the parallelepiped containing the domain and $n$ is the number of
dimensions. We select a regular grid of control points $P$ in the unitary
hypercube and we perturb the space by moving these points. The
displacements, the so called FFD weights, control the basis functions whose
tensor product constitute the deformation map $\hat{T}$. We underline
that it is also possible to move only some points: typically we fix several rows/columns
of control points to obtain desired levels of continuity and to fix certain
parts of the domain.

\begin{figure}[h]
\centering\includegraphics[width=.66\textwidth]{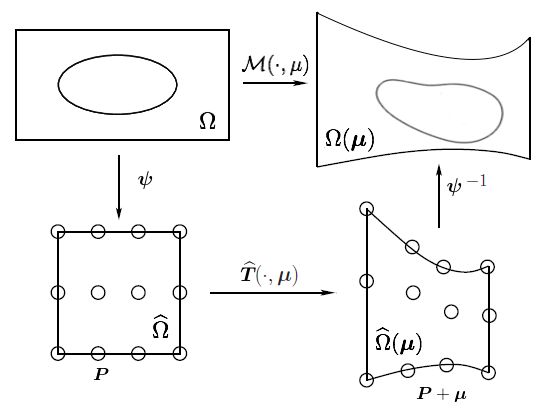}
\caption{Schematic diagram of the free form deformation map
  $\mathcal{M}(\cdot,\mu)$, of the control points $P_{i,k}$, and the resulting
  deformation when applied to the original domain
  $\Omega$. $\mathcal{M}(\cdot,\mu)$ is the composition of the three
  maps presented: $\psi$, $\hat{T}$, and $\psi^{-1}$.}
\label{fig:FFD} 
\end{figure}

Finally, we need the back mapping to the physical domain, that is the map $\psi^{-1}$.
Formally, we obtain the FFD map as the composition of the three maps,
i.e. $\mathcal{M}(\cdot, \mupar) := (\psi \circ \hat{T} \circ
\psi^{-1})(\cdot, \mupar)$, where $\mupar$ refers to the parametric
displacement of the control points (see Figure~\ref{fig:FFD} for a
schematic summary).

It must be remarked that, although FFD is characterized by high
flexibility and easiness of handling, it suffers from some
limitations. The first lies in the fact that the design variables may have
no physical significance: they are defined in a parametric domain that
can not be expressed into a particular unit of measurement by
definition. Moreover all the control points are restricted to lie on a
regular lattice and, in that way, local refinements could not be
performed. 

\begin{figure}[h]
\centering\includegraphics[width=1.\textwidth]{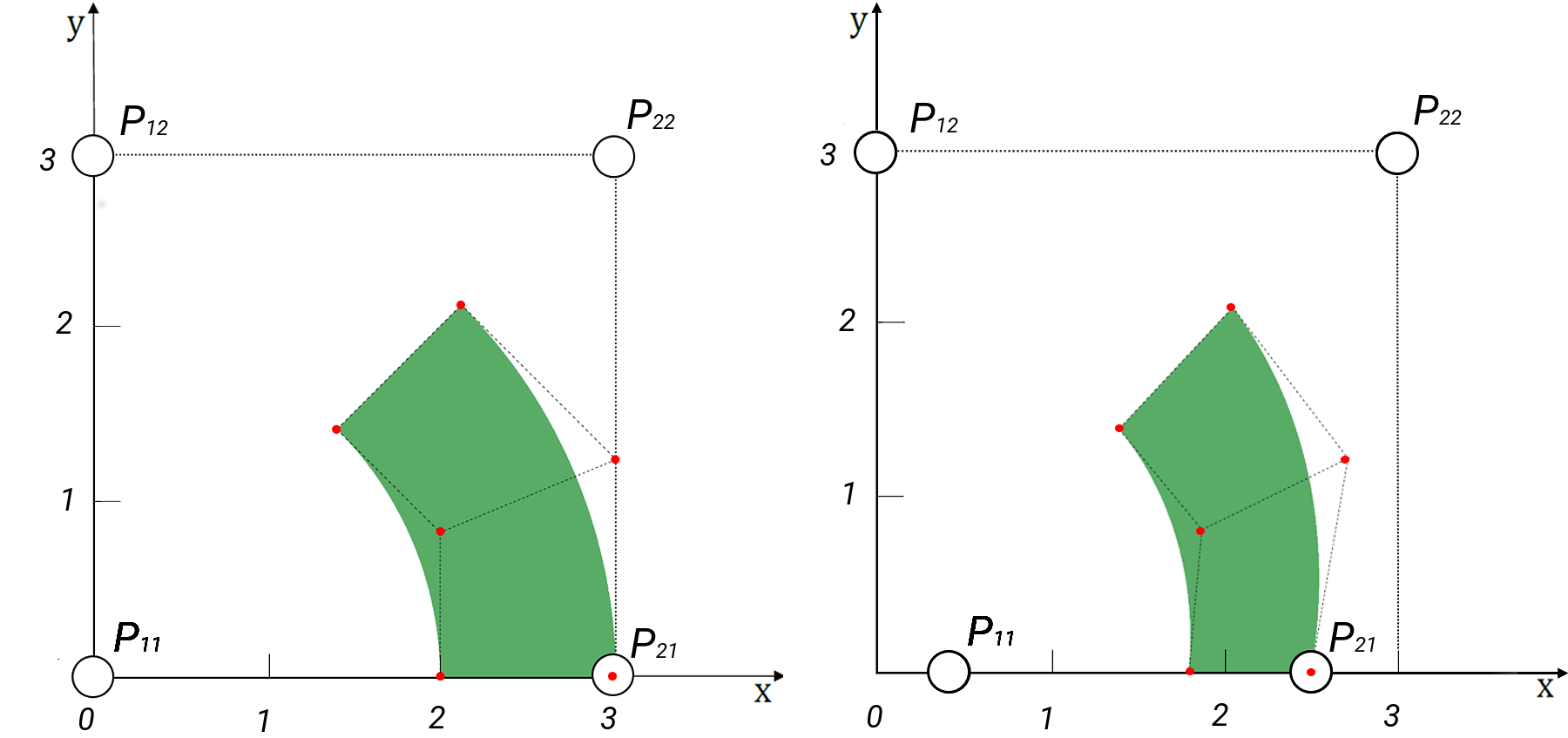}
    \caption{The initial unperturbed domain (left) and an example of deformed
    domain (right) using the FFD technique with $\mupar = (0.4, -0.5)$. The red
    dots are the NURBS control points, the white big dots are the FFD control
    points.}
\label{fig:ffdexample} 
\end{figure}
In this work, we apply the FFD to parametrize the initial 2D domain. We embed
the domain with a square lattice of length 3, using $2 \times 2$ control
points. The lattice origin coincides with the axes origin. We use two different
parameters that are the displacement along the $x$ direction of the
FFD control points $P_{11}$ and $P_{12}$ depicted in
Figure~\ref{fig:ffdexample}. In particular we define $\mathbb{D} :=
[-0.3, 0.3]^2$. We use Bernstein polynomials as basis function to
deform the geometry in the reference domain. In
Figure~\ref{fig:ffdexample} we present on the left the undeformed configuration of
the idealized collector pipe, where in red we highlight the IGA
control points, while the white big dots are the FFD control points. On
the right there is just an example of deformation corresponding to
a displacement of $0.4$ for $P_{11}$ and $-0.5$ for $P_{12}$, for now on
express as $\mupar = (0.4, -0.5)$.

\section{Data-driven reduced order modelling by proper orthogonal
  decomposition with interpolation}
\label{sec:rom}

The reduced basis (RB) method is a computational reduction technique allowing
to quickly and accurately obtain the solution of parametric PDEs. 
The need to solve parameterized differential problems, possibly in a very rapid
calculation time, emerges in various contexts, particularly when we are
interested in characterizing the response of a system in numerous scenarios or
operating
conditions~\cite{rozza2011reduced,rozza2007reduced,peterson1989reduced,quarteroni2014ms}.
The goal of an RB approximation is the representation of the full-order problem
as combination of the (few) essential characteristics of the problem itself. In
this way the dimensions are considerably lower than those of a problem
discretized with a classic Galerkin method.
Any discretization leading to a large system to be solved to achieve a certain
accuracy is referred to as \emph{high fidelity} (or \emph{full order})
\emph{approximation}.
The basic idea of an RB approximation is a computationally efficient solution
of the parametric problem keeping the approximation error lower than a given
tolerance.  In particular, the aim is to approximate the solution of a
parametric PDE using a very small number of degrees of freedom instead of the
large number required by an high fidelity approximation.

To do this, the full order problem is solved only for a few instances of the
input parameters during the computationally expensive \emph{offline} phase. The
so stored snapshots are used in the \emph{online} phase for the approximation
of the solution for any new parameter. The generation of the snapshots
database can be done only once and it is completely decoupled from any new
input-output calculation related to a new parameter. The online phase
exploits the calculations already performed and therefore not
necessary of a large computational power. This 
ensures real-time processing of the problem without having to use high
performance computing infrastructures for
analysis, which can be instead run on a simple laptop with limited
computational power. 

In this work, we adopt a complete data-driven model order reduction called
proper orthogonal decomposition with interpolation (PODI). PODI is applicable
using only the system output --- so also experimental data --- without
requiring the equations of the original problem. Especially in the industrial
context, this is a big benefit: it allows to preserve the \emph{know how} and
to be completely independent from the full-order solver. We list some examples
of PODI
applications~\cite{salmoiraghi2018,bui2003proper,chinesta2016model,hesthaven2016certified}.

PODI aims to approximate the solution manifold by interpolating the snapshots
collected in the offline phase. Since for high-dimensional data the
interpolation can be very expensive, we use proper orthogonal decomposition
(POD) to project the original snapshots onto a low-rank space. POD allows to
define a subspace approximating the original data in an optimal least squares
sense by using the singular value decomposition (SVD)
algorithm~\cite{volkwein2013proper,quarteroni2009numerical,sirovich1987turbulence}.
We consider a set of $n_{\text{train}}$ snapshots ${s_1, \dotsc,
s_{n_{\text{train}}}} = {s(\mu_1), \dotsc ,s(\mu_{n_{\text{train}}})} \in
V^{\mathcal{N}}$, where $V^\mathcal{N}$ is the high-dimensional space and
$\mathcal{N}$ refers to its dimension. We define the snapshots matrix
$\mathbf{S}$ as the matrix that contains the snapshots in the columns
$\mathbf{S} =\begin{bmatrix} s_1 & \dotsc & s_{n_{\text{train}}}\end{bmatrix}$.
    We apply the SVD to $\mathbf{S}$:
\begin{equation}
    \mathbf{S} = \mathbf{V} \boldsymbol{\Sigma} \mathbf{W}^*,
\end{equation}
where
\begin{align}
    \mathbf{V} &= \begin{bmatrix}
        \zeta_1 & \dotsc & \zeta_{n_{\text{train}}}
    \end{bmatrix}\in \mathbb{C}^{\mathcal{N}
    \times n_{\text{train}}},\\
\mathbf{W} &= \begin{bmatrix}
    \Psi_1 & \dotsc & \Psi_{n_{\text{train}}}
    \end{bmatrix}
    \in \mathbb{C}^{{n_{\text{train}}} \times {n_{\text{train}}}},
\end{align}
are orthogonal matrices whose columns are the left and right singular vectors of $\mathbf{S}$ respectively, and
\begin{equation}
    \boldsymbol{\Sigma} = diag(\sigma_1 \dotsc \sigma_{n_{\text{train}}}) \in \mathbb{C}^{{n_{\text{train}}} \times {n_{\text{train}}}},
\end{equation}
is a diagonal matrix such that $\sigma_1 \ge \sigma_2 \ge \dotsc \ge
\sigma_{\text{train}} \ge 0$ are the computed singular values of $\mathbf{S}$.
The POD modes of dimension $N$ are defined as the first $N$ left singular
vectors of $\mathbf{S}$, that correspond to the $N$ largest singular values
\begin{equation}
\mathbf{Z} = \begin{bmatrix}\zeta_1 & \dotsc & \zeta_N
    \end{bmatrix}.
\end{equation} 

Now we project the original snapshots onto the space spanned by the modes: the
snapshots are so described as linear combination of the modes such that
\begin{equation} 
    s_i = \sum_{j=1}^{N} \mathbf{C}_{j,i}\zeta_j 
    \qquad \text{for}\,\,i = 1, \dotsc, n_{\text{train}},
\end{equation}
where the columns of matrix $\mathbf{C}$ are called \emph{modal
coefficients}. We can compute these coefficients as $\mathbf{C} = \mathbf{Z}^T
\mathbf{S}$, where $\mathbf{C} = \begin{bmatrix} c_1 & \dotsc &
c_{n_\text{train}}\end{bmatrix} \in \mathbb{R}^{N \times n_\text{train}}$.
We remark the relation between these coefficients and the parameters; hence we
can interpolate the modal coefficients to compute the coefficient for any new
point belonging to the parameter space. Finally, using the modes, we are able
to approximate the new high-dimensional solution.

Since the PODI technique relies on interpolation, the accuracy of the
approximated solution depends mostly by the chosen interpolation method.

\section{Numerical Results}
\label{sec:results}

In order to construct the reduced order model, we firstly need to sample the
solution manifold using several high-fidelity snapshots.  We select
$100$ different configurations applying the FFD technique to the
initial domain. The parameters $\mupar$ are equispaced in the
parameter space $[-0.3, 0.3] \times
[-0.3, 0.3]$. This strategy allows us to
cover the entire parametric space with a linear interpolation.

\begin{figure}[!h]
\centering
\includegraphics[width=0.45\textwidth]{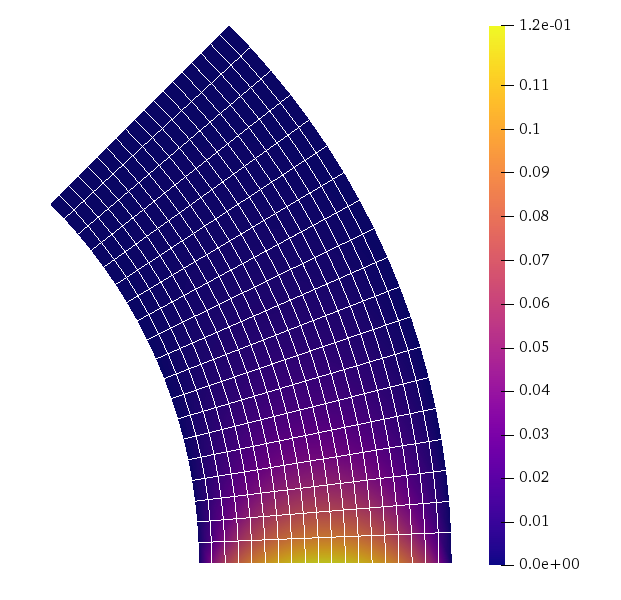}
\caption{The adopted computational grid and the graphical representation of the
    numerical solution of the Laplace problem for the undeformed configuration.}
\label{fig:2D3D}
\end{figure}

For each configuration, IGA is performed testing
\emph{GeoPDEs}~\cite{de2011geopdes,Vazquez_2016}, an open source and free package introduced
in 2010 by Rafael V\'azquez, written in Octave and fully compatible with
Matlab. In GeoPDEs the IGA is efficiently implemented in its classic
Galerkin version. For the resolution of the full-order problem we
created a mesh with 400 degrees of freedom. 
Figure~\ref{fig:2D3D} shows the graphical representation of the numerical
solution.

\begin{figure}[!h]
\centering
\begin{tikzpicture}

\definecolor{color0}{rgb}{0.12156862745098,0.466666666666667,0.705882352941177}

\begin{axis}[
width=.9\textwidth,
height=6cm,
tick align=outside,
tick pos=left,
x grid style={white!69.01960784313725!black},
xlabel={Number of singular values},
xmajorgrids,
xmin=-3.95, xmax=104.95,
y grid style={white!69.01960784313725!black},
ylabel={$^{\sigma}/_{\sigma_1}$},
ymajorgrids,
ymin=1e-10, ymax=6.31094614774382,
ymode=log
]
\addplot [semithick, color0, mark=*, mark size=1, mark options={solid}, forget plot]
table [row sep=\\]{%
1	1 \\
2	0.0266069367530292 \\
3	0.00666716538429062 \\
4	0.00113908028231995 \\
5	0.000335864770088419 \\
6	8.68272115017186e-05 \\
7	4.53996167540263e-05 \\
8	1.419148253871e-05 \\
9	4.49658464627143e-06 \\
10	1.96384699798989e-06 \\
11	8.16735651567433e-07 \\
12	7.63584910874785e-07 \\
13	3.89893143006576e-07 \\
14	3.80748534422903e-07 \\
15	3.6579809519233e-07 \\
16	3.42210618067721e-07 \\
17	3.32984150561027e-07 \\
18	3.10956501310746e-07 \\
19	2.98776137412508e-07 \\
20	2.91318922790681e-07 \\
21	2.77721126415208e-07 \\
22	2.72158122168289e-07 \\
23	2.51047366577065e-07 \\
24	2.43067808664668e-07 \\
25	2.33727597801655e-07 \\
26	2.23877746383615e-07 \\
27	2.19266568318497e-07 \\
28	2.07404269009445e-07 \\
29	2.06638489701005e-07 \\
30	1.97513559015695e-07 \\
31	1.79063048296703e-07 \\
32	1.73372520493746e-07 \\
33	1.59899313002149e-07 \\
34	1.57368732437147e-07 \\
35	1.48643838681325e-07 \\
36	1.42657150001697e-07 \\
37	1.30538726116392e-07 \\
38	1.25627587695671e-07 \\
39	1.19421181387196e-07 \\
40	1.11725896802637e-07 \\
41	1.08016379037939e-07 \\
42	1.01090425140451e-07 \\
43	9.59060006910493e-08 \\
44	9.0780255507734e-08 \\
45	8.5575047341473e-08 \\
46	7.7639211791146e-08 \\
47	7.27088491411874e-08 \\
48	6.61011299588072e-08 \\
49	6.23321177413113e-08 \\
50	5.69927936395169e-08 \\
51	5.33504968799277e-08 \\
52	4.93189839083773e-08 \\
53	4.62187898563934e-08 \\
54	4.2411665887978e-08 \\
55	4.01233304593882e-08 \\
56	3.7950602683685e-08 \\
57	3.66390771374674e-08 \\
58	3.61764693114182e-08 \\
59	3.58796491967991e-08 \\
60	3.38902214696685e-08 \\
61	3.29970745539525e-08 \\
62	3.23807189697768e-08 \\
63	3.1096374583961e-08 \\
64	3.00638165633751e-08 \\
65	2.90471901265689e-08 \\
66	2.82956978513251e-08 \\
67	2.72514897493412e-08 \\
68	2.64379298667225e-08 \\
69	2.54042536226255e-08 \\
70	2.46297199177595e-08 \\
71	2.40931728848093e-08 \\
72	2.27273629921691e-08 \\
73	2.21793360037435e-08 \\
74	2.13904904556511e-08 \\
75	2.07409603165971e-08 \\
76	2.038878547996e-08 \\
77	1.90538658472601e-08 \\
78	1.7597401875838e-08 \\
79	1.73790913301407e-08 \\
80	1.69140042732008e-08 \\
81	1.55497645010532e-08 \\
82	1.46541968352323e-08 \\
83	1.405142917304e-08 \\
84	1.27754023499016e-08 \\
85	1.16290912261983e-08 \\
86	1.04161084623439e-08 \\
87	9.49612772665097e-09 \\
88	9.0421848037042e-09 \\
89	7.58582007240351e-09 \\
90	7.19470766526503e-09 \\
91	6.58948956953866e-09 \\
92	6.30485915447778e-09 \\
93	5.17655258119779e-09 \\
94	4.70072100896797e-09 \\
95	3.82012752557664e-09 \\
};
\end{axis}

\end{tikzpicture}
\caption{The singular values obtained by the snapshots matrix using the SVD
    technique.}
\label{fig:svd}
\end{figure}
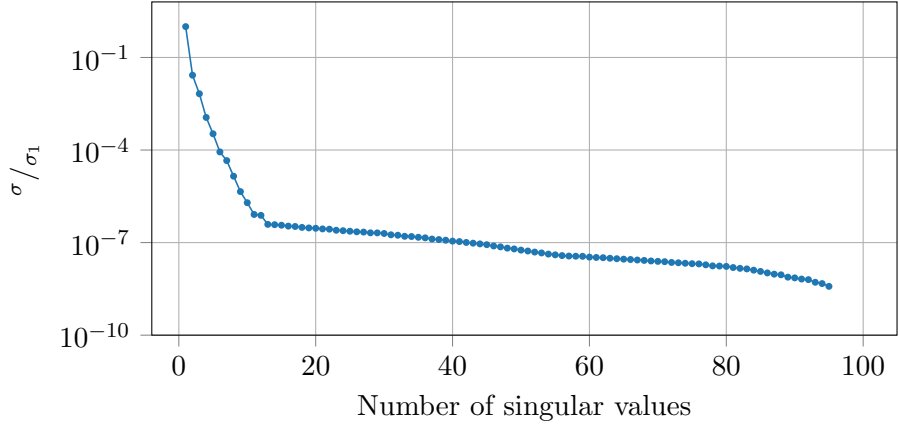

Once the snapshots are collected, we create the reduced order model using the
PODI method. The modes are so computed by applying the SVD algorithm to the
snapshots matrix. We show in Figure~\ref{fig:svd} the obtained singular values:
we note that the first one retains $\sim 96\%$ of the total energy, while
the $10$th singular value is below $10^{-6}$. We expect that even with
only few modes we can generate a reduced order model introducing only a negligible
error.

\begin{figure}[!h]
\centering
\includegraphics[width=.45\textwidth]{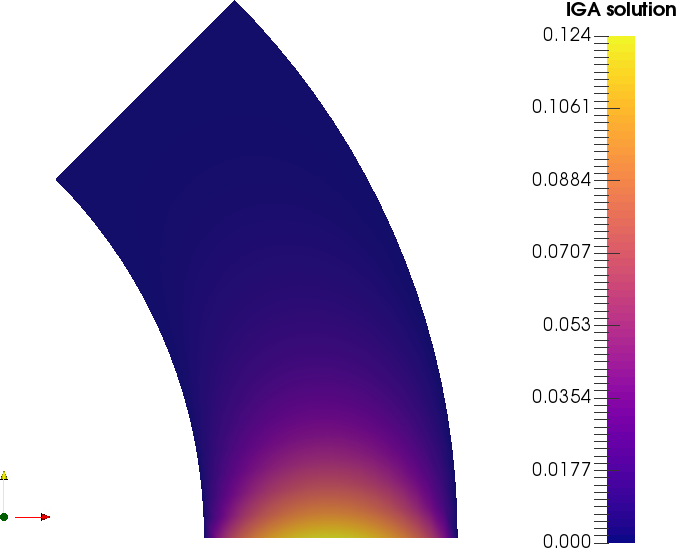}
\includegraphics[width=.45\textwidth]{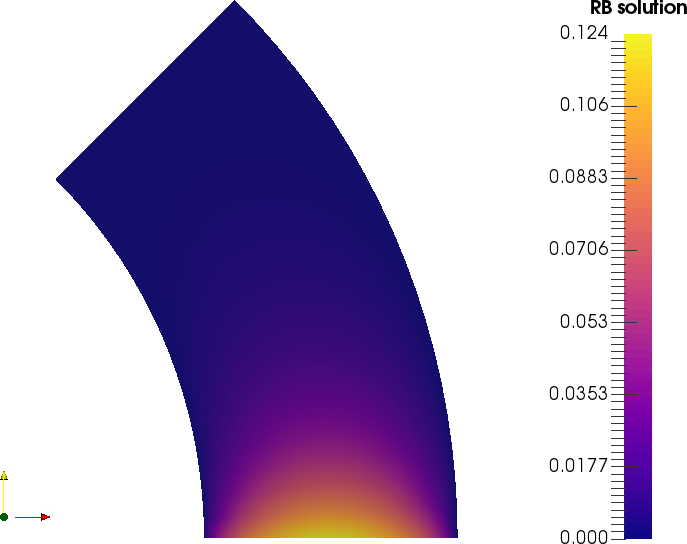}
\caption{Comparison between the full-order solution (left) and the
     reduced order model solution (right) for the undeformed configuration.}
\label{fig:deformation}
\end{figure}

Using the modes, we can calculate the modal coefficients by projecting the
original snapshots. Hence, we can approximate any new solution in the
parametric space trough the interpolation of the modal coefficients. Among the
various interpolation techniques we choose \emph{linear interpolation}. We
report an example where the reduced solution is calculated for the undeformed
object by setting the parameter to zero. In Figure~\ref{fig:deformation} a
visual comparison between the high-fidelity solution and the reduced one 
is presented: it is very intuitive to note that the two solutions are almost
identical. 

\begin{figure}[!h]
\centering
\includegraphics[width=0.5\textwidth]{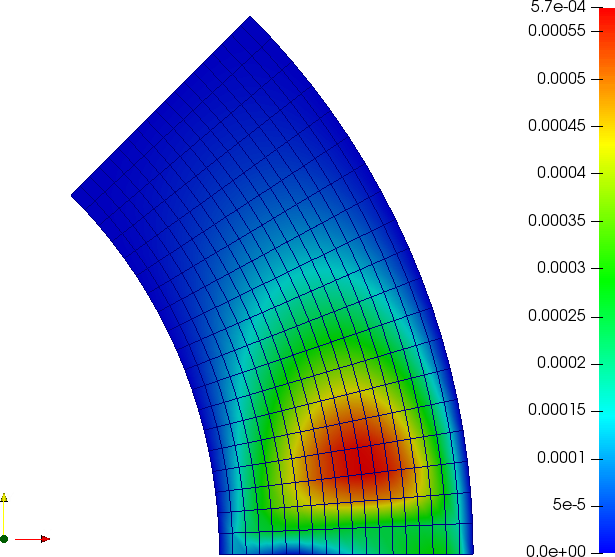}
\caption{Error between the full-order solution and the reduced order solution
    for the unperturbed configuration.
    The different color corresponds to ascending values of the error. The error
    gradually increases from the blue zone to the red zone where there is the
    greatest diffusive heat effect. The results show how, even in the red zone,
    the error assumes acceptable values.}
  \label{fig:differencegraf}
\end{figure}

In Figure~\ref{fig:differencegraf} instead we can see the error 
between the reduced solution and the IGA solution. We calculate the
error $\mathrm{e(\mupar)}$ as follows 
\begin{equation}
\mathrm{e(\mupar)} = |u_{\mathcal{N}}(\mupar) - u_N(\mupar)|.
\end{equation}
The maximum error is around $6 \times 10^{-4}$ so it is
possible to state that it is an acceptable error.

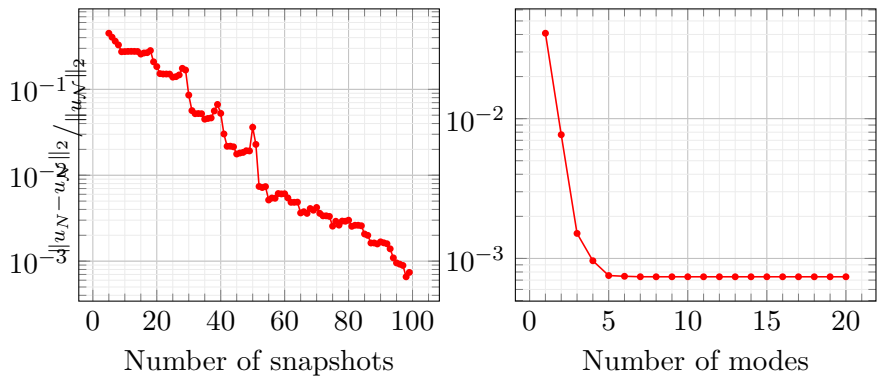
\begin{figure}[!h]
\centering
\begin{tikzpicture}

\begin{groupplot}[group style={group size=2 by 1}]
\nextgroupplot[
width=.5\textwidth,
y label style={at={(axis description cs:0.05,.5)}},
xlabel={Number of snapshots},
ylabel={$^{\|u_N - u_\mathcal{N}\|_2}/_{\|u_\mathcal{N}\|_2}$},
grid=both,
xtick={0,20,...,100},
grid style={line width=.15pt, draw=gray!15},
major grid style={line width=.2pt,draw=gray!50},
minor tick num=3,
xmajorgrids,
ymajorgrids,
ymode=log
]
\addplot [semithick, red, mark=*, mark size=1, mark options={solid}, forget plot]
table [row sep=\\]{%
5	0.449911812990582 \\
6	0.404572636048381 \\
7	0.363565354303628 \\
8	0.328833911448645 \\
9	0.274111658012369 \\
10	0.275303205956197 \\
11	0.276461200673601 \\
12	0.276804858792392 \\
13	0.276600267598185 \\
14	0.275756217345335 \\
15	0.257414453942255 \\
16	0.265091075381544 \\
17	0.26739183371297 \\
18	0.282691324510097 \\
19	0.209120685819603 \\
20	0.184107197763977 \\
21	0.152381300865449 \\
22	0.150772587949068 \\
23	0.15078957318488 \\
24	0.150534762609592 \\
25	0.138767123676949 \\
26	0.140537815414081 \\
27	0.147094942052765 \\
28	0.175638452032969 \\
29	0.167829756827805 \\
30	0.0857648500584303 \\
31	0.0565382642369143 \\
32	0.0520817505827033 \\
33	0.0522421443849487 \\
34	0.0520490636586997 \\
35	0.044938171294647 \\
36	0.045809671672763 \\
37	0.0464785646480516 \\
38	0.0560582921852755 \\
39	0.0667114482249754 \\
40	0.0525602898280551 \\
41	0.0302429760154845 \\
42	0.0217401260999467 \\
43	0.0217912209070218 \\
44	0.0214273454426718 \\
45	0.01767590766636 \\
46	0.0181419879877821 \\
47	0.0184478351419256 \\
48	0.0192854129414914 \\
49	0.0191914949531096 \\
50	0.0361791286482636 \\
51	0.0228587000364364 \\
52	0.00740684797246844 \\
53	0.00719613713757957 \\
54	0.0073731621936017 \\
55	0.00515241725039741 \\
56	0.00542818436882333 \\
57	0.00539450060416036 \\
58	0.00612348925947792 \\
59	0.00606191385308752 \\
60	0.00608474720839682 \\
61	0.00544764611579682 \\
62	0.00483911329382644 \\
63	0.00483911329382652 \\
64	0.00486326806915602 \\
65	0.00364124583017574 \\
66	0.00376178482398796 \\
67	0.00359697997631541 \\
68	0.00409699655719558 \\
69	0.00391588083130021 \\
70	0.00421396093144618 \\
71	0.00359887475882415 \\
72	0.00337011134359315 \\
73	0.0033701113435934 \\
74	0.00331774332761992 \\
75	0.00255591281290352 \\
76	0.00291249159561868 \\
77	0.00262525414084856 \\
78	0.00292679779545679 \\
79	0.00291409019248178 \\
80	0.00299767625067969 \\
81	0.00254457397358195 \\
82	0.00261603140212479 \\
83	0.0026160314021248 \\
84	0.00258165013337853 \\
85	0.00207006955704989 \\
86	0.00199641682673304 \\
87	0.00162640118613441 \\
88	0.00162740975802553 \\
89	0.00158596644644046 \\
90	0.00167725070198238 \\
91	0.00163424984125568 \\
92	0.00159430527783545 \\
93	0.00139138631721693 \\
94	0.00109212976472893 \\
95	0.000953631122502567 \\
96	0.000922158744606224 \\
97	0.000890954003085628 \\
98	0.000657579532211881 \\
99	0.000739259797093588 \\
};
\nextgroupplot[
width=.5\textwidth,
xlabel={Number of modes},
grid=both,
xtick={0,5,...,20},
grid style={line width=.15pt, draw=gray!15},
major grid style={line width=.2pt,draw=gray!50},
minor tick num=4,
xmajorgrids,
ymajorgrids,
ymode=log
]
\addplot [semithick, red, mark=*, mark size=1, mark options={solid}, forget plot]
table [row sep=\\]{%
1	0.0408250625062223 \\
2	0.00767376866346784 \\
3	0.00151060212772228 \\
4	0.000962589774885566 \\
5	0.000754629013773978 \\
6	0.000745402923892498 \\
7	0.000739801704457556 \\
8	0.000739350534893891 \\
9	0.000739271179618598 \\
10	0.0007392594288499 \\
11	0.000739255570025353 \\
12	0.000739254355370234 \\
13	0.000739254667146183 \\
14	0.000739255034314937 \\
15	0.000739254954930179 \\
16	0.000739255097296065 \\
17	0.000739255481547788 \\
18	0.000739256820539752 \\
19	0.00073925731982177 \\
20	0.00073925785981931 \\
};
\end{groupplot}

\end{tikzpicture}
\caption{A posteriori $L^2$ relative error between the reduced order
  solution and the high fidelity one, computed on the test dataset composed by
  20 uniformly distributed random samples in the parameter space. On
  the left the error with respect to the dimension of the offline
  database. On the right the error trend varying the number of modes
  selected. With only 4 modes we obtain an average relative error
  below $10^{-3}$.}
\label{fig:errortrend}
\end{figure}

We can also evaluate the a posteriori error committed on a test dataset.  A posteriori error
estimation allows to minimize the dimension $N$ of the snapshot database used
to generate the reduced space and to quantify the error of the approximation
with respect to the number of modes selected. The error is calculated computing
the relative $L^2$ norm of the difference between the approximated solution
obtained using PODI approach and the IGA truth solution, on a
test dataset composed by high fidelity solutions corresponding to 20
uniformly distributed random samples in the parameter space. 
The plots in Figure~\ref{fig:errortrend} show the relative error
against the dimension of the database and the number of modes. We see
that 100 samples and only 4 modes are enough for an average error below
$10^{-3}$. We refer to~\cite{DevaudRozza2017} for a posteriori error bounds in
an RB-IGA setting.

Finally, we can evaluate the performance improvement
obtained using ROM by calculating the speedup $S_p$ as 
\begin{equation}
S_p = \frac{u_{\mathcal{N}}(s)}{u_N(s)},
\end{equation}
where we divide the time in seconds needed to compute the full order solution
by the time needed for the reduced one. Due to the different size of the
systems, the difference of computational time is remarkable even if, for this
test-case, the full-order problem is very simple. We measured the computational
time required by the two techniques on the same machine, and for different
parameter values, and we obtained a mean speedup of approximately $1000$.
Concerning the software involved, for the model order reduction we adopted
EZyRB~\cite{demo18ezyrb}, which is a Python library for ROM, based on baricentric
triangulation for the selection of the parameter points and on POD for the
selection of the modes. The software uses a non-intrusive approach in which the
solutions are projected on the low dimensional space then interpolated for the
approximation of the solution.

\section{Conclusions and future developments}
\label{sec:the_end}

In this work we presented a complete non-intrusive computational pipeline
involving
geometrical parameterization through free form deformation,
isogeometric analysis, and reduced order model, for fast and reliable
field evaluation. We applied this pipeline to a diffusion problem in
an idealized 2D collector pipe. We used a data-driven non-intrusive
approach for the model order reduction, that is the proper orthogonal
decomposition with interpolation. This setting, even if tested on a
simple problem, will allow us to deal with more complex industrial CAD
files, since we used a geometrical parameterization technique
independent from the object of interest, and the ROM chosen uses only
the snapshots of the IGA high fidelity simulations.

Results and speedup achieved look promising to continue with the
implementation of more complex problems on 3D geometries. The
effectiveness of an RB approach would be exploited even better
increasing the complexity of the simulation in cases where a large
number of analysis has to be computed, e.g. in parameter optimization
studies. The developed RB-IGA method is thus interesting from both
academic and industrial points of view. As a matter of fact, since IGA
is directly interfaced with CAD, an undergoing development of the
work is the implementation of a dedicated software based on the RB-IGA
method, allowing real-time evaluations of  outputs of interest for
different NURBS parameterizations.

\section*{Acknowledgments}
This work was partially funded by the project SOPHYA,
``Seakeeping Of Planing Hull YAchts'', supported by Regione
FVG, POR-FESR 2014-2020, Piano Operativo Regionale Fondo Europeo per
lo Sviluppo Regionale, and partially supported by European Union Funding for
Research and Innovation --- Horizon 2020 Program --- in the framework
of European Research Council Executive Agency: H2020 ERC CoG 2015
AROMA-CFD project 681447 ``Advanced Reduced Order Methods with
Applications in Computational Fluid Dynamics'' P.I. Gianluigi Rozza. 

This work was also partially supported by Fondazione Cariplo - Regione
Lombardia through the project ``Verso nuovi strumenti di simulazione
super veloci ed accurati basati sull'analisi isogeometrica'', within
the program RST - rafforzamento. 


\end{document}